\input amstex


\nopagenumbers
\frenchspacing
\font\erm=cmr10
\font\mrm=cmr10 scaled\magstep1
\hfuzz= 200pt
\vfuzz= 1pt
\vsize = 8.6 true in

\headline{\ifnum\pageno=1\hfil\else
{\ifodd\pageno\rightheadline \else \leftheadline\fi}\fi}
\def\rightheadline{\hfil\hskip 2.3in\erm 
q-zeros \hskip 2.3in
\erm\folio}
\def\leftheadline{\erm\folio \hskip 2.3in 
\erm I.Cherednik\hskip 2.3in\hfil }

\mrm
\baselineskip=12pt

\centerline{ Calculating  zeros of a $q$-zeta function numerically}
\vskip 0.3cm
\centerline{Ivan Cherednik}
\vskip 0.5cm

This note is an appendix to the previous paper
``On q-analogues of Riemann's zeta''.
We will give an output of the computer program 
evaluating the zeros of the following modification
of the ``plus-sharp'' $q$-zeta function:
$$
\eqalignno{
&\zeta^\#(k)\ =\ 
\sum_{j=0}^\infty\,\prod_{l=1}^j
{(1-e^{-(l+2k-1)/a})(1-e^{(l+k)/a})\over
 (1-e^{-(l+k-1)/a})(1-e^{l/a})}\,
{\exp({dk^2\over 4a})+1\over \exp({d(k+j)^2\over 4a})+1}.
&(1)
}
$$
Here $q=\exp(-1/a), \ a,d>0.$

These zeros are  the deformations of the zeros
$k=yi$ of the classical $\zeta(k+1/2).$
They are considered in the horizontal strip  
$$
\eqalignno{
&K(a,d)=\{k,\ 0<\Im k< 2\epsilon \}\hbox{\ for\ }
\epsilon =\sqrt{{\pi a\over 2d}}.
&(2)
}
$$
We use the formula for the linear approximation $z_a$ of
the zero $z$ corresponding to the classical $yi.$
$$
\eqalignno{
& z_a= yi(1-{ 4d^{-1}(1/2+yi)\zeta_+(3/2+yi)-d(-1+yi)\zeta_+(-1/2+yi)
\over 12 a\zeta_+'(1/2+yi)})\cr
&\hbox{\ for\ } \zeta_+(s)=(1-2^{1-s})\zeta(s),\ \zeta_+'(s)=d \zeta_+(s)/ds.
&(3)
}
$$

Given a classical zero, the program carefully verifies the existence
of the zero $z$ of $\zeta^\#$ in a neighborhood of $z_a$ 
and evaluates $z$ numerically.
We did not try to reach  high accuracy when calculating $z.$ 
The procedure is  a sequence of integrations
over  diminishing  rectangles around the consecutive approximations.

The change of the 
$\arg(\zeta^\#)$ gives the number of the
zeros (poles) inside the considered  rectangle. In the tables below,
this number is the difference between the first and the last 
{\it angles} divided by $2\pi.$

We note that  some (but not all) zeros can be found by means of the Newton
method or its modifications. In the output, 
by {\it iterations} we mean the Newton iterations
applied to the zero found by means of the integrations.

Actually the program is universal. It can be used for any function
analytic in the upper half-plane if the  approximations of the
zeros are known, especially 
if the latter are near the imaginary axis. 
For instance, it can be applied in the soliton theory.
However its main purpose is theoretical.
\vskip 0.3cm

\vfil
Notation:

1) $y*I$ = zero of $\zeta,\ z$ = zero of $\zeta^\#, 
\ za$ = its linear approximation,

2) $zn$ = the center of the rectangle of $re\times im-$
dimensions $2*rd\times 2*rad,$ 

3) $zna$ is either $za$ or the value of $z$ after the last {\it good}
integration,

4) $char$ is $1$ minus the number of zeros in the considered rectangle,

5) $fo=\max\{1+[2(|\alpha-\alpha'|-1)]\}$ if $|\alpha-\alpha'|>1$ for all
pairs of the consecutive {\it angles},

6) $vv$ is $|\zeta^\#(z)|$ divided by  $|\zeta^\#(zn)|$,

6a) in the final list of zeros, it is divided by the value at $za$,

7) $10*de$ is a certain upper boundary  of the error for $z,$

8) $[\,b*\sqrt{a/d}\,]$ is the number of terms in the sum (1), 

9) $4*c$ is the number of points on the rectangle used for the integration, 

10) {\it good}  mainly means that $char=0, fo\le 2,  vv<0.8,$ 

10a) the resulting $z$ is inside the integration rectangle, and

10b) $|\zeta^\#(z)|$ is not bigger than the previous {\it good}
    absolute value,

11) {\it very good} is the second consecutive {\it good} with $fo\le 1$
    and admissible  $de,\ vv.$ 

\vskip 0.2cm

The program is entirely automated. 
It changes the center $zn$ and the dimensions
$2rd\times 2rad$ after each integration. If the result of the 
integration is {\it good} then $rd$ and $rad$ are divided by $2.$
A {\it very good} integration concludes  the consideration of a
given classical zero. The program adds more points to the rectangles,
more integrations, and more variants if necessary. We note that
each {\it angle} in the output is actually the sum of $4$ angles 
over  $4$ points located on the different sides of a given
integration rectangle. 

Let $a=750,\ d=2,$ initial $c=4.$ The $b$ are calculated automatically.
\vskip 0.5cm
\vfil

PROGRAM: Q-ZEROS of SHARP ZETA  
\vskip 0.2cm 
\vskip 0.2cm 

CLASSICAL ZEROS AND APPROXIMATIONS:

d= 2  a= 750  ALL ZEROS TILL 48.5406:
\vskip 0.2cm

1  y= 14.1347  za= 0.1303 + 14.1465 I  b= 15

2  y= 21.0220  za= 0.3504 + 21.0771 I  b= 15

3  y= 25.0109  za= 0.5745 + 24.9643 I  b= 15

4  y= 30.4249  za= 0.9134 + 30.4077 I  b= 15

5  y= 32.9351  za= 1.0998 + 33.0854 I  b= 15

6  y= 37.5862  za= 1.7675 + 38.1895 I  b= 20

7  y= 40.9187  za= 1.9141 + 40.7816 I  b= 20

8  y= 43.3271  za= 2.4497 + 43.3138 I  b= 20

9  y= 48.0052  za= 3.1103 + 47.5578 I  b= 20
\vskip 0.2cm 
\vskip 0.3cm
\vfil

FINAL LIST OF Q-ZEROS:
\vskip 0.2cm
 
very good 1  14.1347  z: 0.1304 + 14.1450 I
                                                          
  za= 0.1303 + 14.1465 I  de= 0.000001   vv= 0.000001
\vskip 0.2cm  
very good 2  21.0220  z: 0.3514 + 21.0702 I
                                                        
  za= 0.3504 + 21.0771 I  de= 0.000001   vv= 0.000001
\vskip 0.2cm  
very good 3  25.0109  z: 0.5641 + 24.9586 I
                                       
  za= 0.5745 + 24.9643 I  de= 0.000001   vv= 0.000001
\vskip 0.2cm  
very good 4  30.4249  z: 0.9046 + 30.4014 I
                                          
  za= 0.9134 + 30.4077 I  de= 0.000001   vv= 0.000019
\vskip 0.2cm  
very good 5  32.9351  z: 1.1051 + 33.0341 I
                                        
  za= 1.0998 + 33.0854 I   de= 0.000001  vv= 0.000116
\vskip 0.2cm  
very good 6  37.5862  z: 1.6449 + 37.9659 I
 
  za= 1.7675 + 38.1895 I  de= 0.000040   vv= 0.001053
\vskip 0.2cm  
very good 7  40.9187  z: 1.9080 + 40.8119 I
 
  za= 1.9141 + 40.7816 I  de= 0.000031   vv= 0.013871
\vskip 0.2cm  
very good 8  43.3271  z: 2.2860 + 43.2485 I
 
  za= 2.4497 + 43.3138 I  de= 0.000030   vv= 0.001586
\vskip 0.2cm  
very good 9  48.0052  z: 2.9259 + 47.8424 I
  
  za= 3.1103 + 47.5578 I  de= 0.000005   vv= 0.003049

\vfill\eject

\erm
\baselineskip=12pt

    \obeylines
    \nopagenumbers
    \frenchspacing
    \parindent=0pt
    \let\lr=L \newbox\leftcolumn
    \hsize=2.9 true in
    \newdimen\fullhsize
    \fullhsize=6.2 true in
    \def\fullline{\hbox to\fullhsize}
    \def\columnbox{\leftline{\pagebody}}
    \output={\if L\lr
       \global\setbox\leftcolumn=\columnbox \global\let\lr=R
       \else \doubleformat \global\let\lr=L\fi
      \ifnum\outputpenalty>-20000 \else\dosupereject\fi}
    \def\doubleformat{\shipout\vbox{\makeheadline
        \fullline{\box\leftcolumn\hfil\columnbox}
        \makefootline}
      \advancepageno}
    \def\endcol{\par\vfill\eject}
    \def\endpage{\endcol\supereject
      \if R\lr\null\vfill\eject\fi}

\comment
1  y= 14.1347  za= 0.130263 + 14.1465 I  nzb= 15
2  y= 21.022  za= 0.35044 + 21.0771 I  nzb= 15
3  y= 25.0109  za= 0.574478 + 24.9643 I  nzb= 15
4  y= 30.4249  za= 0.913386 + 30.4077 I  nzb= 15
5  y= 32.9351  za= 1.09982 + 33.0854 I  nzb= 15
6  y= 37.5862  za= 1.76753 + 38.1895 I  nzb= 20
7  y= 40.9187  za= 1.9141 + 40.7816 I  nzb= 20
8  y= 43.3271  za= 2.44967 + 43.3138 I  nzb= 20
9  y= 48.0052  za= 3.11028 + 47.5578 I  nzb= 20
\endcomment
 
VARIANT= 1  c= 4
\vskip 0.3cm 
 
\vskip 0.3cm  no= 1  y= 14.1347 b= 15 c= 4
zna= 0.130263 + 14.1465 I  
zn= 0.130263 + 14.1465 I
rd= 0.0477578 rad= 0.0238789
angles over the rd*rad rectangle: 
1     -2.79498
1  1  -3.00821
1  2  -3.25288
1  3  -3.53768
2     -3.87423
2  1  -4.27666
2  2  -4.75796
2  3  -5.31955
3     -5.93574
3  1  -6.55215
3  2  -7.11419
3  3  -7.59589
4     -7.99857
4  1  -8.33527
4  2  -8.62015
4  3  -8.86488
5     -9.07817
char= 0. fo= 0
rd= 0.0477578 rad= 0.0238789
z= 0.130488 + 14.1452 I 
zn= 0.130263 + 14.1465 I
vv= 0.136129
 
good

\vskip 0.3cm  no= 1  y= 14.1347 b= 15 c= 4
zna= 0.130488 + 14.1452 I  
zn= 0.130461 + 14.1453 I
rd= 0.0238789 rad= 0.0119394
angles over the rd*rad rectangle: 
1     2.26652
1  1  2.05311
1  2  1.80819
1  3  1.52302
2     1.18598
2  1  0.782973
2  2  0.301227
2  3  -0.26021
3     -0.875348
3  1  -1.49038
3  2  -2.0516
3  3  -2.53316
4     -2.93607
4  1  -3.27309
4  2  -3.55826
4  3  -3.80322
5     -4.01667
char= 0. fo= 0
rd= 0.0238789 rad= 0.0119394
z= 0.130414 + 14.145 I 
zn= 0.130461 + 14.1453 I
vv= 0.250379
 
good
 
\vskip 0.3cm second try:

\vskip 0.3cm  no= 1  y= 14.1347 b= 15 c= 6
zna= 0.130414 + 14.145 I  
zn= 0.130414 + 14.145 I
rd= 0.0119394 rad= 0.00596972
angles over the rd*rad rectangle: 
1     1.05302
2     0.409533
3     -0.582258
3  1  -0.912566
3  2  -1.2788
3  3  -1.67502
4     -2.08862
4  1  -2.50221
4  2  -2.89841
4  3  -3.26462
5     -3.5949
6     -4.58666
7     -5.23017
char= 0. fo= 0
rd= 0.0119394 rad= 0.00596972
z= 0.130402 + 14.145 I 
zn= 0.130414 + 14.145 I
vv= 0.495182
 
good

\vskip 0.3cm  no= 1  y= 14.1347 b= 15 c= 6
zna= 0.130402 + 14.145 I  
zn= 0.130406 + 14.145 I
rd= 0.00596972 rad= 0.00298486
angles over the rd*rad rectangle: 
1     1.02567
2     0.382187
3     -0.609613
3  1  -0.939931
3  2  -1.30618
3  3  -1.70242
4     -2.11602
4  1  -2.52962
4  2  -2.92581
4  3  -3.29201
5     -3.62228
6     -4.61402
7     -5.25752
char= 0. fo= 0
rd= 0.00596972 rad= 0.00298486
z= 0.130392 + 14.145 I 
zn= 0.130406 + 14.145 I
vv= 0.250002
 
good
very good
iterations:  
z=0.130389 + 14.145 I 
vv=0.0000188769

\vskip 0.3cm  no= 2  y= 21.022 b= 15 c= 4
zna= 0.35044 + 21.0771 I  
zn= 0.35044 + 21.0771 I
rd= 0.129534 rad= 0.064767
angles over the rd*rad rectangle: 
1     -2.56518
1  1  -2.77792
1  2  -3.02189
1  3  -3.30573
2     -3.64099
2  1  -4.04187
2  2  -4.52194
2  3  -5.08381
3     -5.70266
3  1  -6.32267
3  2  -6.88694
3  3  -7.36905
4     -7.77118
4  1  -8.10708
4  2  -8.39126
4  3  -8.63545
5     -8.84837
char= 0. fo= 0
rd= 0.129534 rad= 0.064767
z= 0.352258 + 21.0712 I 
zn= 0.35044 + 21.0771 I
vv= 0.177848
 
good

\vskip 0.3cm  no= 2  y= 21.022 b= 15 c= 4
zna= 0.352258 + 21.0712 I  
zn= 0.351984 + 21.0721 I
rd= 0.064767 rad= 0.0323835
angles over the rd*rad rectangle: 
1     2.01358
1  1  1.80038
1  2  1.5557
1  3  1.2708
2     0.933995
2  1  0.531067
2  2  0.0490028
2  3  -0.513352
3     -1.1297
3  1  -1.74542
3  2  -2.30648
3  3  -2.78749
4     -3.18985
4  1  -3.52646
4  2  -3.81139
4  3  -4.05621
5     -4.26961
char= 0. fo= 0
rd= 0.064767 rad= 0.0323835
z= 0.35165 + 21.0705 I 
zn= 0.351984 + 21.0721 I
vv= 0.250897
 
good
 
\vskip 0.3cm second try:

\vskip 0.3cm  no= 2  y= 21.022 b= 15 c= 6
zna= 0.35165 + 21.0705 I  
zn= 0.35165 + 21.0705 I
rd= 0.0323835 rad= 0.0161918
angles over the rd*rad rectangle: 
1     0.334116
2     -0.309329
3     -1.3012
3  1  -1.63157
3  2  -1.99788
3  3  -2.39417
4     -2.8078
4  1  -3.22137
4  2  -3.61751
4  3  -3.98365
5     -4.31387
6     -5.30554
7     -5.94907
char= 0. fo= 0
rd= 0.0323835 rad= 0.0161918
z= 0.351549 + 21.0703 I 
zn= 0.35165 + 21.0705 I
vv= 0.496271
 
good

\vskip 0.3cm  no= 2  y= 21.022 b= 15 c= 6
zna= 0.351549 + 21.0703 I  
zn= 0.351583 + 21.0704 I
rd= 0.0161918 rad= 0.00809588
angles over the rd*rad rectangle: 
1     0.30616
2     -0.337286
3     -1.32924
3  1  -1.65967
3  2  -2.02606
3  3  -2.42241
4     -2.83607
4  1  -3.24963
4  2  -3.64572
4  3  -4.0118
5     -4.34196
6     -5.33351
7     -5.97703
char= 0. fo= 0
rd= 0.0161918 rad= 0.00809588
z= 0.351473 + 21.0702 I 
zn= 0.351583 + 21.0704 I
vv= 0.249984
 
good
very good
iterations:  
z=0.351447 + 21.0702 I 
vv=0.00027944

\vskip 0.3cm  no= 3  y= 25.0109 b= 15 c= 4
zna= 0.574478 + 24.9643 I  
zn= 0.574478 + 24.9643 I
rd= 0.210457 rad= 0.105228
angles over the rd*rad rectangle: 
1     -1.13776
1  1  -1.35014
1  2  -1.59447
1  3  -1.8797
2     -2.21773
2  1  -2.62288
2  2  -3.10757
2  3  -3.671
3     -4.28463
3  1  -4.89496
3  2  -5.45168
3  3  -5.93109
4     -6.33392
4  1  -6.67202
4  2  -6.95879
4  3  -7.20557
5     -7.42095
char= 0. fo= 0
rd= 0.210457 rad= 0.105228
z= 0.565775 + 24.9608 I 
zn= 0.574478 + 24.9643 I
vv= 0.226358
 
good

\vskip 0.3cm  no= 3  y= 25.0109 b= 15 c= 4
zna= 0.565775 + 24.9608 I  
zn= 0.567381 + 24.9614 I
rd= 0.105228 rad= 0.0526142
angles over the rd*rad rectangle: 
1     -0.289756
1  1  -0.502799
1  2  -0.747503
1  3  -1.03271
2     -1.37022
2  1  -1.77434
2  2  -2.2579
2  3  -2.82116
3     -3.43662
3  1  -4.04992
3  2  -4.60876
3  3  -5.08871
4     -5.49097
4  1  -5.82798
4  2  -6.11349
4  3  -6.35893
5     -6.57294
char= 0. fo= 0
rd= 0.105228 rad= 0.0526142
z= 0.564553 + 24.9592 I 
zn= 0.567381 + 24.9614 I
vv= 0.252223
 
good
 
\vskip 0.3cm second try:

\vskip 0.3cm  no= 3  y= 25.0109 b= 15 c= 6
zna= 0.564553 + 24.9592 I  
zn= 0.564553 + 24.9592 I
rd= 0.0526142 rad= 0.0263071
angles over the rd*rad rectangle: 
1     0.445348
2     -0.197914
3     -1.18976
3  1  -1.52017
3  2  -1.88655
3  3  -2.28291
4     -2.69659
4  1  -3.11016
4  2  -3.50627
4  3  -3.87236
5     -4.20253
6     -5.19416
7     -5.83784
char= 0. fo= 0
rd= 0.0526142 rad= 0.0263071
z= 0.564342 + 24.9589 I 
zn= 0.564553 + 24.9592 I
vv= 0.498128
 
good

\vskip 0.3cm  no= 3  y= 25.0109 b= 15 c= 6
zna= 0.564342 + 24.9589 I  
zn= 0.564412 + 24.959 I
rd= 0.0263071 rad= 0.0131536
angles over the rd*rad rectangle: 
1     0.405824
2     -0.237536
3     -1.22959
3  1  -1.56014
3  2  -1.92665
3  3  -2.32313
4     -2.73687
4  1  -3.15043
4  2  -3.54645
4  3  -3.91242
5     -4.24247
6     -5.23382
7     -5.87736
char= 0. fo= 0
rd= 0.0263071 rad= 0.0131536
z= 0.564177 + 24.9587 I 
zn= 0.564412 + 24.959 I
vv= 0.250036
 
good
very good
iterations:  
z=0.564122 + 24.9586 I 
vv=0.000542266

\vskip 0.3cm  no= 4  y= 30.4249 b= 15 c= 4
zna= 0.913386 + 30.4077 I  
zn= 0.913386 + 30.4077 I
rd= 0.333581 rad= 0.16679
angles over the rd*rad rectangle: 
1     -0.681479
1  1  -0.893326
1  2  -1.13696
1  3  -1.42118
2     -1.75773
2  1  -2.16079
2  2  -2.64308
2  3  -3.20501
3     -3.81972
3  1  -4.43334
3  2  -4.99326
3  3  -5.47442
4     -5.87778
4  1  -6.21587
4  2  -6.5025
4  3  -6.74921
5     -6.96466
char= 0. fo= 0
rd= 0.333581 rad= 0.16679
z= 0.909368 + 30.4063 I 
zn= 0.913386 + 30.4077 I
vv= 0.627134
 
good

\vskip 0.3cm  no= 4  y= 30.4249 b= 15 c= 4
zna= 0.909368 + 30.4063 I  
zn= 0.910917 + 30.4068 I
rd= 0.16679 rad= 0.0833952
angles over the rd*rad rectangle: 
1     -0.314662
1  1  -0.527545
1  2  -0.772143
1  3  -1.05735
2     -1.39509
2  1  -1.79975
2  2  -2.28424
2  3  -2.84845
3     -3.46416
3  1  -4.07668
3  2  -4.63436
3  3  -5.11349
4     -5.51536
4  1  -5.85231
4  2  -6.13793
4  3  -6.38359
5     -6.59785
char= 0. fo= 0
rd= 0.16679 rad= 0.0833952
z= 0.905809 + 30.4026 I 
zn= 0.910917 + 30.4068 I
vv= 0.244024
 
good
 
\vskip 0.3cm second try:

\vskip 0.3cm  no= 4  y= 30.4249 b= 15 c= 6
zna= 0.905809 + 30.4026 I  
zn= 0.905809 + 30.4026 I
rd= 0.0833952 rad= 0.0416976
angles over the rd*rad rectangle: 
1     -0.101447
2     -0.744697
3     -1.73688
3  1  -2.06747
3  2  -2.43404
3  3  -2.83055
4     -3.2442
4  1  -3.65765
4  2  -4.05353
4  3  -4.41934
5     -4.74937
6     -5.7408
7     -6.38463
char= 0. fo= 0
rd= 0.0833952 rad= 0.0416976
z= 0.905214 + 30.402 I 
zn= 0.905809 + 30.4026 I
vv= 0.497814
 
good

\vskip 0.3cm  no= 4  y= 30.4249 b= 15 c= 6
zna= 0.905214 + 30.402 I  
zn= 0.905412 + 30.4022 I
rd= 0.0416976 rad= 0.0208488
angles over the rd*rad rectangle: 
1     -0.151109
2     -0.794509
3     -1.78724
3  1  -2.11812
3  2  -2.48499
3  3  -2.88165
4     -3.29543
4  1  -3.7088
4  2  -4.10433
4  3  -4.46984
5     -4.79955
6     -5.79056
7     -6.43429
char= 0. fo= 0
rd= 0.0416976 rad= 0.0208488
z= 0.904746 + 30.4015 I 
zn= 0.905412 + 30.4022 I
vv= 0.249506
 
good

\vskip 0.3cm  no= 5  y= 32.9351 b= 15 c= 4
zna= 1.09982 + 33.0854 I  
zn= 1.09982 + 33.0854 I
rd= 0.405333 rad= 0.202666
angles over the rd*rad rectangle: 
1     -2.76933
1  1  -2.97572
1  2  -3.21238
1  3  -3.48756
2     -3.81255
2  1  -4.20243
2  2  -4.67463
2  3  -5.24022
3     -5.87965
3  1  -6.52496
3  2  -7.10238
3  3  -7.58516
4     -7.98297
4  1  -8.31438
4  2  -8.59568
4  3  -8.83886
5     -9.05251
char= 0. fo= 0
rd= 0.405333 rad= 0.202666
z= 1.11278 + 33.0399 I 
zn= 1.09982 + 33.0854 I
vv= 0.176273
 
good

\vskip 0.3cm  no= 5  y= 32.9351 b= 15 c= 4
zna= 1.11278 + 33.0399 I  
zn= 1.11083 + 33.0467 I
rd= 0.202666 rad= 0.101333
angles over the rd*rad rectangle: 
1     1.44051
1  1  1.22913
1  2  0.986248
1  3  0.702958
2     0.36725
2  1  -0.0359012
2  2  -0.520682
2  3  -1.08877
3     -1.71112
3  1  -2.32861
3  2  -2.88712
3  3  -3.36441
4     -3.76398
4  1  -4.09905
4  2  -4.38348
4  3  -4.62856
5     -4.84267
char= 0. fo= 0
rd= 0.202666 rad= 0.101333
z= 1.10703 + 33.0356 I 
zn= 1.11083 + 33.0467 I
vv= 0.252102
 
good
 
\vskip 0.3cm second try:

\vskip 0.3cm  no= 5  y= 32.9351 b= 15 c= 6
zna= 1.10703 + 33.0356 I  
zn= 1.10703 + 33.0356 I
rd= 0.101333 rad= 0.0506666
angles over the rd*rad rectangle: 
1     -0.417857
2     -1.06077
3     -2.05313
3  1  -2.38381
3  2  -2.75047
3  3  -3.14699
4     -3.56074
4  1  -3.97409
4  2  -4.36961
4  3  -4.73524
5     -5.06523
6     -6.0567
7     -6.70105
char= 0. fo= 0
rd= 0.101333 rad= 0.0506666
z= 1.10609 + 33.0348 I 
zn= 1.10703 + 33.0356 I
vv= 0.502104
 
good

\vskip 0.3cm  no= 5  y= 32.9351 b= 15 c= 6
zna= 1.10609 + 33.0348 I  
zn= 1.1064 + 33.035 I
rd= 0.0506666 rad= 0.0253333
angles over the rd*rad rectangle: 
1     -0.459909
2     -1.10327
3     -2.0966
3  1  -2.42765
3  2  -2.79487
3  3  -3.19158
4     -3.6054
4  1  -4.01845
4  2  -4.41347
4  3  -4.77897
5     -5.1083
6     -6.09915
7     -6.74312
char= 0. fo= 0
rd= 0.0506666 rad= 0.0253333
z= 1.10535 + 33.0342 I 
zn= 1.1064 + 33.035 I
vv= 0.250262
 
good
very good
iterations:  
z=1.1051 + 33.034 I 
vv=0.0207832

\vskip 0.3cm  no= 6  y= 37.5862 b= 20 c= 4
zna= 1.76753 + 38.1895 I  
zn= 1.76753 + 38.1895 I
rd= 0.5 rad= 0.25
angles over the rd*rad rectangle: 
1     1.31082
2     0.586877
2  1  0.22623
2  2  -1.07829
2 new 0.586877
2  1  0.435554
2  2  0.22623
2  3  -0.149968
2  4  -1.07829
2  5  -2.09925
2  6  -2.52703
2  7  -2.76634
3     -2.9474
3  1  -3.25021
3  2  -3.52701
3  3  -3.793
4     -4.05083
5     -4.97236
char= 0. fo= 1
rd= 0.5 rad= 0.25
z= 1.6622 + 37.9732 I 
zn= 1.76753 + 38.1895 I
vv= 0.0459034
 
good

\vskip 0.3cm  no= 6  y= 37.5862 b= 20 c= 4
zna= 1.69731 + 38.0453 I  
zn= 1.66683 + 37.9827 I
rd= 0.25 rad= 0.125
angles over the rd*rad rectangle: 
1     -0.525067
1  1  -0.73675
1  2  -0.981341
1  3  -1.26771
2     -1.60898
2  1  -2.02142
2  2  -2.51782
2  3  -3.09305
3     -3.70924
3  1  -4.31009
3  2  -4.85343
3  3  -5.32232
4     -5.72005
4  1  -6.05667
4  2  -6.34362
4  3  -6.59129
5     -6.80825
char= 0. fo= 0
rd= 0.25 rad= 0.125
z= 1.64767 + 37.9678 I 
zn= 1.66683 + 37.9827 I
vv= 0.192703
 
good
 
\vskip 0.3cm second try:

\vskip 0.3cm  no= 6  y= 37.5862 b= 20 c= 6
zna= 1.64767 + 37.9678 I  
zn= 1.64767 + 37.9678 I
rd= 0.125 rad= 0.0625
angles over the rd*rad rectangle: 
1     -0.723796
2     -1.36698
3     -2.35965
3  1  -2.69088
3  2  -3.05835
3  3  -3.4537
4     -3.86757
4  1  -4.28133
4  2  -4.67449
4  3  -5.04026
5     -5.37102
6     -6.36244
7     -7.00698
char= 0. fo= 0
rd= 0.125 rad= 0.0625
z= 1.6462 + 37.9666 I 
zn= 1.64767 + 37.9678 I
vv= 0.512907
 
good

\vskip 0.3cm  no= 6  y= 37.5862 b= 20 c= 6
zna= 1.6462 + 37.9666 I  
zn= 1.6467 + 37.967 I
rd= 0.0625 rad= 0.03125
angles over the rd*rad rectangle: 
1     -0.861441
2     -1.50514
3     -2.5004
3  1  -2.83125
3  2  -3.1998
3  3  -3.59624
4     -4.01044
4  1  -4.42139
4  2  -4.81491
4  3  -5.18145
5     -5.50947
6     -6.49934
7     -7.14463
char= 0. fo= 0
rd= 0.0625 rad= 0.03125
z= 1.64494 + 37.9659 I 
zn= 1.6467 + 37.967 I
vv= 0.232099
 
good
very good
iterations do not work

\vskip 0.3cm  no= 7  y= 40.9187 b= 20 c= 4
zna= 1.9141 + 40.7816 I  
zn= 1.9141 + 40.7816 I
rd= 0.5 rad= 0.25
angles over the rd*rad rectangle: 
1     -2.35413
1  1  -2.56065
1  2  -2.79839
1  3  -3.07707
2     -3.40794
2  1  -3.80464
2  2  -4.28063
2  3  -4.84018
3     -5.459
3  1  -6.08126
3  2  -6.64756
3  3  -7.13198
4     -7.53793
4  1  -7.87772
4  2  -8.16663
4  3  -8.41707
5     -8.63731
char= 0. fo= 0
rd= 0.5 rad= 0.25
z= 1.92012 + 40.8209 I 
zn= 1.9141 + 40.7816 I
vv= 0.544958
 
good

\vskip 0.3cm  no= 7  y= 40.9187 b= 20 c= 4
zna= 1.92012 + 40.8209 I  
zn= 1.918 + 40.8071 I
rd= 0.25 rad= 0.125
angles over the rd*rad rectangle: 
1     1.45608
1  1  1.24379
1  2  0.999723
1  3  0.714207
2     0.37875
2  1  -0.02443
2  2  -0.502504
2  3  -1.06163
3     -1.67244
3  1  -2.28494
3  2  -2.849
3  3  -3.33314
4     -3.73864
4  1  -4.07763
4  2  -4.36536
4  3  -4.61165
5     -4.82714
char= 0. fo= 0
rd= 0.25 rad= 0.125
z= 1.91055 + 40.8136 I 
zn= 1.918 + 40.8071 I
vv= 0.225127
 
good
 
\vskip 0.3cm second try:

\vskip 0.3cm  no= 7  y= 40.9187 b= 20 c= 6
zna= 1.91055 + 40.8136 I  
zn= 1.91055 + 40.8136 I
rd= 0.125 rad= 0.0625
angles over the rd*rad rectangle: 
1     -0.744738
2     -1.38803
3     -2.37935
3  1  -2.71058
3  2  -3.078
3  3  -3.47346
4     -3.88863
4  1  -4.30088
4  2  -4.6961
4  3  -5.06075
5     -5.39018
6     -6.38198
7     -7.02792
char= 0. fo= 0
rd= 0.125 rad= 0.0625
z= 1.90905 + 40.8125 I 
zn= 1.91055 + 40.8136 I
vv= 0.463925
 
good

\vskip 0.3cm  no= 7  y= 40.9187 b= 20 c= 6
zna= 1.90905 + 40.8125 I  
zn= 1.90953 + 40.8129 I
rd= 0.0625 rad= 0.03125
angles over the rd*rad rectangle: 
1     -0.921727
2     -1.56266
3     -2.55779
3  1  -2.89042
3  2  -3.26241
3  3  -3.65551
4     -4.0676
4  1  -4.48351
4  2  -4.87042
4  3  -5.24271
5     -5.57032
6     -6.55652
7     -7.20491
char= 0. fo= 0
rd= 0.0625 rad= 0.03125
z= 1.908 + 40.8119 I 
zn= 1.90953 + 40.8129 I
vv= 0.243719
 
good
very good
iterations do not work

\vskip 0.3cm  no= 8  y= 43.3271 b= 20 c= 4
zna= 2.44967 + 43.3138 I  
zn= 2.44967 + 43.3138 I
rd= 0.5 rad= 0.25
angles over the rd*rad rectangle: 
1     -1.61602
1  1  -1.83468
1  2  -2.10201
1  3  -2.43246
2     -2.8472
2  1  -3.35542
2  2  -3.9102
2  3  -4.4511
3     -4.9509
3  1  -5.41813
3  2  -5.86163
3  3  -6.28312
4     -6.6786
4  1  -7.03872
4  2  -7.36158
4  3  -7.64693
5     -7.89921
char= 0. fo= 0
rd= 0.5 rad= 0.25
z= 2.29315 + 43.2565 I 
zn= 2.44967 + 43.3138 I
vv= 0.049899
 
good

\vskip 0.3cm  no= 8  y= 43.3271 b= 20 c= 4
zna= 2.29315 + 43.2565 I  
zn= 2.30059 + 43.2592 I
rd= 0.25 rad= 0.125
angles over the rd*rad rectangle: 
1     -0.632114
1  1  -0.844933
1  2  -1.08782
1  3  -1.37466
2     -1.7122
2  1  -2.11869
2  2  -2.60296
2  3  -3.16815
3     -3.7831
3  1  -4.39779
3  2  -4.95012
3  3  -5.42485
4     -5.82729
4  1  -6.16506
4  2  -6.45158
4  3  -6.70026
5     -6.9153
char= 0. fo= 0
rd= 0.25 rad= 0.125
z= 2.28858 + 43.2501 I 
zn= 2.30059 + 43.2592 I
vv= 0.305191
 
good
 
\vskip 0.3cm second try:

\vskip 0.3cm  no= 8  y= 43.3271 b= 20 c= 6
zna= 2.28858 + 43.2501 I  
zn= 2.28858 + 43.2501 I
rd= 0.125 rad= 0.0625
angles over the rd*rad rectangle: 
1     -1.00531
2     -1.64511
3     -2.63818
3  1  -2.96883
3  2  -3.33458
3  3  -3.72935
4     -4.15314
4  1  -4.56094
4  2  -4.95928
4  3  -5.32061
5     -5.6524
6     -6.641
7     -7.2885
char= 0. fo= 0
rd= 0.125 rad= 0.0625
z= 2.28718 + 43.249 I 
zn= 2.28858 + 43.2501 I
vv= 0.467616
 
good

\vskip 0.3cm  no= 8  y= 43.3271 b= 20 c= 6
zna= 2.28718 + 43.249 I  
zn= 2.28762 + 43.2494 I
rd= 0.0625 rad= 0.03125
angles over the rd*rad rectangle: 
1     -1.20849
2     -1.858
3     -2.84843
3  1  -3.18015
3  2  -3.55535
3  3  -3.94177
4     -4.35293
4  1  -4.77126
4  2  -5.16281
4  3  -5.53312
5     -5.85965
6     -6.85195
7     -7.49168
char= 0. fo= 0
rd= 0.0625 rad= 0.03125
z= 2.28604 + 43.2485 I 
zn= 2.28762 + 43.2494 I
vv= 0.22265
 
good
very good
iterations do not work

\vskip 0.3cm  no= 9  y= 48.0052 b= 20 c= 4
zna= 3.11028 + 47.5578 I  
zn= 3.11028 + 47.5578 I
rd= 0.5 rad= 0.25
angles over the rd*rad rectangle: 
1     -1.38493
2     -2.21348
3     -2.9527
3  1  -3.05309
3  2  -2.96753
3  3  -2.01601
4     -1.18241
5     -1.38493
char= 1. fo= 0
rd= 0.5 rad= 0.25
z= 2.98924 + 47.3299 I 
zn= 3.11028 + 47.5578 I
vv= 0.895232

\vskip 0.3cm  no= 9  y= 48.0052 b= 20 c= 4
zna= 3.11028 + 47.5578 I  
zn= 3.08002 + 47.5009 I
rd= 1. rad= 0.5
angles over the rd*rad rectangle: 
1     -1.91196
2     -2.75572
2  1  -3.02683
2  2  -3.33321
2  3  -3.70806
3     -4.21951
3  1  -5.04389
3  2  -6.01395
3  3  -6.6674
4     -7.11581
4  1  -7.46582
4  2  -7.7524
4  3  -7.98942
5     -8.19515
char= 0. fo= 0
rd= 1. rad= 0.5
z= 3.28887 + 47.8094 I 
zn= 3.08002 + 47.5009 I
vv= 2.16026
 
\vskip 0.3cm second try:

\vskip 0.3cm  no= 9  y= 48.0052 b= 20 c= 6
zna= 3.11028 + 47.5578 I  
zn= 3.11028 + 47.5578 I
rd= 0.812049 rad= 0.406024
angles over the rd*rad rectangle: 
1     -1.28591
2     -1.83358
3     -2.53892
4     -3.50803
4  1  -3.88532
4  2  -4.4221
4  3  -5.11249
5     -5.72593
5  1  -6.15537
5  2  -6.47188
5  3  -6.71562
6     -6.92932
7     -7.56797
char= 0.000180379 fo= 0
rd= 0.812049 rad= 0.406024
z= 3.05445 + 47.9338 I 
zn= 3.11028 + 47.5578 I
vv= 0.843354

\vskip 0.3cm  no= 9  y= 48.0052 b= 20 c= 6
zna= 3.11028 + 47.5578 I  
zn= 3.08468 + 47.7302 I
rd= 0.609037 rad= 0.304518
angles over the rd*rad rectangle: 
1     0.371705
2     -0.26843
3     -1.14232
3  1  -1.40267
3  2  -1.68278
3  3  -1.97787
4     -2.30382
4  1  -2.66226
4  2  -3.05396
4  3  -3.47895
5     -3.88863
5  1  -4.28275
5  2  -4.61402
5  3  -4.90167
6     -5.15322
7     -5.91129
char= 0. fo= 0
rd= 0.609037 rad= 0.304518
z= 2.96001 + 47.8662 I 
zn= 3.08468 + 47.7302 I
vv= 0.225039
 
good
 
\vskip 0.3cm  
VARIANT= 2  c= 6
\vskip 0.3cm 
 
\vskip 0.3cm  no= 4  y= 30.4249 b= 15 c= 6
zna= 0.904746 + 30.4015 I  
zn= 0.904746 + 30.4015 I
rd= 0.00251499 rad= 0.00125749
angles over the rd*rad rectangle: 
1     -0.214344
2     -0.8561
2  1  -1.06313
2  2  -1.29739
2  3  -1.5615
3     -1.85971
3  1  -2.19599
3  2  -2.56697
3  3  -2.96901
4     -3.38244
4  1  -3.79389
4  2  -4.18433
4  3  -4.54339
5     -4.86709
6     -5.85195
7     -6.49735
char= 0. fo= 0
rd= 0.00251499 rad= 0.00125749
z= 0.90459 + 30.4014 I 
zn= 0.904746 + 30.4015 I
vv= 0.00484321
 
good

\vskip 0.3cm  no= 4  y= 30.4249 b= 15 c= 6
zna= 0.90459 + 30.4014 I  
zn= 0.904591 + 30.4014 I
rd= 0.00125749 rad= 0.000628747
angles over the rd*rad rectangle: 
1     -2.11878
2     -2.76408
3     -3.75334
3  1  -4.0849
3  2  -4.45175
3  3  -4.84679
4     -5.25976
4  1  -5.67336
4  2  -6.06929
4  3  -6.43586
5     -6.76651
6     -7.75825
7     -8.40196
char= 0. fo= 0
rd= 0.00125749 rad= 0.000628747
z= 0.90459 + 30.4014 I 
zn= 0.904591 + 30.4014 I
vv= 0.206104
 
good
very good
iterations do not work

\vskip 0.3cm  no= 9  y= 48.0052 b= 20 c= 6
zna= 2.96001 + 47.8662 I  
zn= 2.96001 + 47.8662 I
rd= 0.573699 rad= 0.286849
angles over the rd*rad rectangle: 
1     -1.03279
2     -1.57682
3     -2.48568
3  1  -2.79493
3  2  -3.14634
3  3  -3.5278
4     -3.92849
4  1  -4.33505
4  2  -4.71914
4  3  -5.08949
5     -5.42839
5  1  -5.73976
5  2  -6.02718
5  3  -6.29151
6     -6.5358
7     -7.31639
char= 0. fo= 0
rd= 0.573699 rad= 0.286849
z= 2.95364 + 47.8503 I 
zn= 2.96001 + 47.8662 I
vv= 0.660047
 
good

\vskip 0.3cm  no= 9  y= 48.0052 b= 20 c= 6
zna= 2.95364 + 47.8503 I  
zn= 2.95617 + 47.8566 I
rd= 0.286849 rad= 0.143425
angles over the rd*rad rectangle: 
1     -1.35265
2     -1.97964
3     -2.97403
3  1  -3.31052
3  2  -3.68183
3  3  -4.07544
4     -4.48391
4  1  -4.8804
4  2  -5.25409
4  3  -5.62115
5     -5.94686
5  1  -6.23352
5  2  -6.49998
5  3  -6.74385
6     -6.95322
7     -7.63583
char= 0. fo= 0
rd= 0.286849 rad= 0.143425
z= 2.93263 + 47.8449 I 
zn= 2.95617 + 47.8566 I
vv= 0.246155
 
good

\vskip 0.3cm  no= 9  y= 48.0052 b= 20 c= 6
zna= 2.93263 + 47.8449 I  
zn= 2.93728 + 47.8472 I
rd= 0.143425 rad= 0.0717123
angles over the rd*rad rectangle: 
1     -1.42112
2     -2.05974
3     -3.04324
3  1  -3.39571
3  2  -3.75441
3  3  -4.16734
4     -4.56347
4  1  -4.96311
4  2  -5.35354
4  3  -5.72723
5     -6.07143
6     -7.04948
7     -7.70311
char= 0.000189447 fo= 0
rd= 0.143425 rad= 0.0717123
z= 2.92799 + 47.8432 I 
zn= 2.93728 + 47.8472 I
vv= 0.383844
 
good
 
\vskip 0.3cm second try:

\vskip 0.3cm  no= 9  y= 48.0052 b= 20 c= 9
zna= 2.92799 + 47.8432 I  
zn= 2.92799 + 47.8432 I
rd= 0.0717123 rad= 0.0358562
angles over the rd*rad rectangle: 
1     -0.426895
2     -0.849362
3     -1.34599
4     -2.0717
5     -3.01717
5  1  -3.25357
5  2  -3.57341
5  3  -3.8698
6     -4.10193
7     -5.0508
8     -5.76695
9     -6.29659
10     -6.70618
char= 0.000621347 fo= 0
rd= 0.0717123 rad= 0.0358562
z= 2.92617 + 47.8424 I 
zn= 2.92799 + 47.8432 I
vv= 0.528196
 
good

\vskip 0.3cm  no= 9  y= 48.0052 b= 20 c= 9
zna= 2.92617 + 47.8424 I  
zn= 2.9268 + 47.8427 I
rd= 0.0358562 rad= 0.0179281
angles over the rd*rad rectangle: 
1     -1.0289
2     -1.42916
3     -1.97926
4     -2.72265
5     -3.64737
5  1  -3.90746
5  2  -4.19868
5  3  -4.47038
6     -4.76644
6  1  -4.95629
6  2  -5.30599
6  3  -5.48123
7     -5.79693
8     -6.38006
9     -6.92474
10     -7.31209
char= 0. fo= 0
rd= 0.0358562 rad= 0.0179281
z= 2.92592 + 47.8424 I 
zn= 2.9268 + 47.8427 I
vv= 0.487768
 
good
very good
iterations do not work
 
\comment 
ZEROS:
 
very good 1  14.1347 z: 0.130389 + 14.145 I6
                                       -6             -8
za= 0.130263 + 14.1465 I de= 1.00002 10   vv= 7.965 10
 
very good 2  21.022 z: 0.351447 + 21.0702 I6
                                      -6              -6
za= 0.35044 + 21.0771 I de= 1.00028 10   vv= 1.5469 10
 
very good 3  25.0109 z: 0.564122 + 24.9586 I6
                                       -6               -6
za= 0.574478 + 24.9643 I de= 1.00301 10   vv= 3.85598 10
 
very good 4  30.4249 z: 0.90459 + 30.4014 I6
                                       -6
za= 0.913386 + 30.4077 I de= 1.00873 10   vv= 0.0000189741
 
very good 5  32.9351 z: 1.1051 + 33.034 I6
                                      -6
za= 1.09982 + 33.0854 I de= 1.30273 10   vv= 0.000116055
 
very good 6  37.5862 z: 1.64494 + 37.9659 I6
za= 1.76753 + 38.1895 I de= 0.0000395443 vv= 0.00105304
 
very good 7  40.9187 z: 1.908 + 40.8119 I6
za= 1.9141 + 40.7816 I de= 0.0000314597 vv= 0.0138717
 
very good 8  43.3271 z: 2.28604 + 43.2485 I6
za= 2.44967 + 43.3138 I de= 0.0000302622 vv= 0.00158553
 
very good 9  48.0052 z: 2.92592 + 47.8424 I6
                                      -6
za= 3.11028 + 47.5578 I de= 6.18724 10   vv= 0.0030494
\endcomment

\bye
\end